\documentclass[10pt,twoside,reqno,centertags,english]{amsart}
\usepackage{hyperref}
\usepackage[active]{srcltx}

\usepackage{amsmath,amssymb}
\usepackage[T1]{fontenc}
\usepackage{color}
\usepackage{textcomp}

\newtheorem*{theo*}{Theorem}{\bfseries}{\itshape}

\newcommand{\eps}{\varepsilon}
\newcommand{\R}{\mathbb{R}}
\newcommand{\W}{\mathcal{W}}
\newcommand{\D}{\mathcal{D}}
\newcommand{\N}{\mathcal{N}}
\newcommand{\B}{\mathfrak{B}}
\newcommand{\E}{\mathcal{E}}
\newcommand{\FR}{\mathsf{FR}}
\renewcommand{\S}{\mathcal{S}}
\renewcommand{\P}{\mathcal{P}}
\newcommand{\Pp}{\mathbb{P}}
\newcommand{\Stat}{\mathfrak{S}}

\newcommand{\rd}{\mathrm{d}}
\newcommand{\grad}{\operatorname{grad}}
\DeclareMathOperator{\dive}{div}
\DeclareMathOperator{\tr}{tr}
\newcommand{\Id}{\operatorname{Id}}

\newcommand{\qtext}[1]{\quad\mbox{#1}\quad}
\newcommand{\qqtext}[1]{\qquad\mbox{#1}\qquad}
\begin{document}
\title{Schr\"odinger encounters Fisher and Rao: a survey}
%
%
\author[L.~Monsaingeon]{L\'eonard Monsaingeon}
\address[L.~Monsaingeon]{GFM Universidade de Lisboa, Campo Grande, Edif\'icio C6, 1749-016 Lisboa, Portugal
\& IECL Universit\'e de Lorraine, F-54506 Vandoeuvre-l\`es-Nancy Cedex, FRANCE
}
\email{leonard.monsaingeon@univ-lorraine.fr}

\author[D.~Vorotnikov]{Dmitry Vorotnikov}
\address[D.~Vorotnikov]{University of Coimbra, CMUC, Department of
Mathematics, 3001-501 Coimbra, Portugal}{}
\email{mitvorot@mat.uc.pt}
\begin{abstract}
In this short note we review the dynamical Schr\"odinger problem on the non-commutative Fisher-Rao space of positive semi-definite matrix-valued measures.
The presentation is meant to be self-contained, and we discuss in particular connections with Gaussian optimal transport, entropy, and quantum Fisher information.
\end{abstract}
\maketitle              
%
%

The abstract Schr\"odinger problem on a generic Riemannian manifold $(M,g)$ can be formulated as the dynamical minimization problem
\begin{equation}
\min\limits_{q}\Bigg\{\int_0^1 |\dot q_t|^2 \rd t + \eps^2\int_0^1 |\nabla V(q_t)|^2 \rd t
\qqtext{s.t. }
q|_{t=0,1}=q_0,q_1
\Bigg\}.
\label{eq:geometric_schro}
\end{equation}
The unknown curves $q=(q_t)_{t\in[0,1]}$ take values in $M$ and interpolate between the prescribed endpoints $q_0,q_1\in M$, the potential $V:M\to \R$ is given, and $\eps>0$ is a regularization parameter.
Clearly when $\eps\to 0$ this $\eps$-problem is expected to converge in some sense to the geodesic problem on $M$.

The original Schr\"odinger problem \cite{schro}, or rather, its dynamical equivalent (sometimes referred to as the Yasue problem \cite{Leger,chen_survey}), can be rewritten as a particular instanciation of \eqref{eq:geometric_schro} when $(M,g)$ is the Wasserstein space of probability measures, the potential $V$ is given by the Boltzmann entropy $H(\rho)=\int \rho\log\rho$, and its Wasserstein gradient $I(\rho)=\|\nabla_\mathcal W H(\rho)\|^2=\int |\nabla\log\rho|^2\rho$ is the Fisher information functional.
This was only recently recognized as an \emph{entropic regularization} of the Monge-Kantorovich problem \cite{leonard} and led to spectacular developments in computational optimal transport \cite{cuturi_peyre}.
Recently, attempts have been made to develop an optimal transport theory for \emph{quantum} objects, namely (symmetric positive-definite) matrix-valued measures, see \cite{chen} and references therein.
In \cite{MV} we studied \eqref{eq:geometric_schro} precisely in the corresponding noncommutative Fisher-Rao space of matrix-measures.
In this note we aim at providing a comprehensive introduction to this specific setting, in particular we wish to emphasize the construction based on Gaussian optimal transport that will lead to some specific entropy and quantum Fisher information functionals later on.

Section~\ref{sec:OT_schro} briefly reviews classical optimal transport and the original Schr\"odinger problem.
In Section~\ref{sec:Bures-Wass} we discuss Gaussian optimal transport and the corresponding Bures-Wasserstein geometry.
This then suggests the natural construction of the Fisher-Rao space detailed in Section~\ref{sec:FR}.
We introduce in Section~\ref{sec:extended_entropy} the resulting entropy and Fisher information functionals dictated by the previous construction.
We also compute explicitly the induced heat flow, and finally discuss the Fisher-Rao-Schr\"odinger problem.

\section{Optimal transport and the Schr\"odinger problem}
\label{sec:OT_schro}
In its fluid-mechanical \emph{Benamou-Brenier} formulation \cite{BB}, the quadratic \emph{Wasserstein distance} between probability measures $\rho_0,\rho_1\in\P(\Omega)$ over a smooth domain $\Omega\subset \R^d$ reads
\begin{multline}
\W^2(\rho_0,\rho_1)=\min\limits_{\rho,v}\Bigg\{
\int_0^1 \int_{\Omega} |v_t(x)|^2\rho_t(x)\rd x \,\rd t
\\
\mbox{s.t.}\quad
\partial_t\rho_t + \dive(\rho_t v_t)=0
\mbox{ in }(0,1)\times \Omega
\qtext{and}\rho|_{t=0,1}=\rho_{0,1}
\Bigg\}
\label{eq:def_W2_dynamic}
\end{multline}
(supplemented with homogeneous no-flux boundary conditions on $\partial\Omega$ if needed in order to ensure mass conservation).
This can be seen as an optimal control problem, where the control $v=v_t(x)$ is used to drive the system from $\rho_0$ to $\rho_1$ while minimizing the overall kinetic energy.
We refer to \cite{villani_small} for a gentle yet comprehensive introduction to optimal transport.


As originally formulated by E. Schr\"odinger himself in \cite{schro}, the Schr\"odinger problem roughly consists in determining the most likely evolution of a system for $t\in[0,1]$ in an ambient noisy environment at temperature $\eps>0$, given the observation of its statistical distribution at times $t=0$ and $t=1$.
Here we shall rather focus on the equivalent dynamical problem
\begin{multline}
\min\limits_{\rho,v}\Bigg\{
\int_0^1 \int_\Omega |v_t(x)|^2 \rho_t(x)\rd x \,\rd t+
\eps^2\int_0^1\int_\Omega |\nabla\log\rho_t(x)|^2\rho_t(x)\rd x \,\rd t
\\
\mbox{s.t.}\quad
\partial_t\rho_t + \dive(\rho_t v_t)=0
\mbox{ in }(0,1)\times \Omega
\qtext{and}\rho|_{t=0,1}=\rho_{0,1}
\Bigg\}
\label{eq:def_Weps_dynamic}
\end{multline}
(again supplemented with Neumann boundary conditions on $\partial\Omega$ if needed).
This is sometimes called the Yasue problem, cf. \cite{Leger,chen_survey}.
It is known \cite{Mikami,leonard,BM} that the noisy problem Gamma-converges towards deterministic optimal transport in the small-temperature limit $\eps\to 0$, and this actually holds in a much more general metric setting than just optimal transport \cite{MVT}.
We also refer to \cite{zambrini,EQM} for connections with Euclidean quantum mechanics, and to the survey \cite{chen_survey} for an optimal-control perspective.

Following F.Otto \cite{otto_geometry}, an important feature of optimal transport is that one can view $\P(\Omega)$ as a (formal) Riemannian manifold, whose Riemannian distance coincides with the Wasserstein distance \eqref{eq:def_W2_dynamic}.
This relies upon the identification of infinitesimal variations $\xi$ with (gradients of) Kantorovich potentials $\phi:\Omega\to \R$ by uniquely selecting the velocity field $v=\nabla\phi$ with minimal kinetic energy $\int_\Omega \rho|v|^2$, given $-\dive(\rho v)=\xi$.
Practically speaking, this means that tangent vectors $\xi\in T_\rho\P$ at a point $\rho$ are identified with scalar potentials $\phi_\xi$ via the Onsager operator $\Delta_\rho=\dive(\rho\nabla \cdot)$ and the elliptic equation
\begin{equation}
-\dive(\rho\nabla\phi_\xi)=\xi,
\label{eq:Onsager_W}
\end{equation}
see \cite{villani_small} for details.
Accordingly, the linear heat flow $\partial_t\rho=\Delta\rho$ can be identified \cite{JKO} as the Wasserstein gradient flow $\frac{d\rho}{dt}=-\grad_\W H(\rho)$ of the \emph{Boltzmann entropy}
$$
H(\rho)=\int_\Omega\rho(x)\log(\rho(x))\rd x.
$$
Here $\grad_\W$ denotes  the gradient computed with respect to Otto's Riemannian structure.
Moreover, the \emph{Fisher information} functional
\begin{equation}
F(\rho)=\int_\Omega |\nabla\log\rho(x)|^2\rho(x)\, \rd x=\|\grad_\W H(\rho)\|^2
\label{eq:def_Fisher_Wasserstein}
\end{equation}
appearing in \eqref{eq:def_Weps_dynamic} coincides with the squared gradient of the entropy, or equivalently with the dissipation rate $-\frac{dH}{dt}=F$ of $H$ along its own gradient flow.
This identification \eqref{eq:def_Fisher_Wasserstein} shows that \eqref{eq:def_Weps_dynamic} can indeed be written as a particular case of \eqref{eq:geometric_schro} with respect to the Wasserstein geometry.

\section{The Bures-Wasserstein distance and Gaussian optimal transport}
\label{sec:Bures-Wass}
When $\rho_0=\N(m_0,A_0),\rho_1=\N(m_1,A_1)$ are Gaussian measures with means $m_i\in\R^d$ and covariance matrices $A_i\in \S^{++}_d(\R)$ (the space of symmetric positive definite matrices), the optimal transport problem \eqref{eq:def_W2_dynamic} is explicitly solvable \cite{takatsu2011wasserstein} and the Wasserstein distance can be computed as
\begin{equation}
\W^2(\rho_0,\rho_1)=|m_1-m_0|^2 + \B^2(A_0,A_1).
\label{eq:W2_Bures_Gaussian_OT}
\end{equation}
Here $\B$ is the \emph{Bures distance}
\begin{equation}
\B^2(A_0,A_1)=\min\limits_{R R^t=\Id } |A_1^{\frac 12}-RA_0^{\frac 12}|^2=\tr A_0 + \tr A_1 -2\tr\left(  \left(A_0^{\frac 12}A_1A_0^{\frac 12}\right)^{\frac 12}\right),
\label{eq:def_Bures_static}
\end{equation}
and $|M|^2=\tr(M M^t)$ is the Euclidean norm of a matrix $M\in \R^{d\times d}$ corresponding to the Frobenius scalar product $\langle M,N\rangle=\tr(MN^t)$.
The Bures distance, sometimes also called the \emph{Helstrom metric}, can be considered as a quantum generalization of the usual Fisher information metric \cite{bures}.
Since the Euclidean geometry of the translational part $|m_1-m_0|^2$ in \eqref{eq:W2_Bures_Gaussian_OT} is trivially flat, we restrict for simplicity to centered Gaussians $m_i=0$.
We denote accordingly
$$
\Stat_0=\Bigg\{
\rho\in\mathcal P(\R^d):\qquad \rho=\N(0,A) \mbox{ for some }A\in\S^{++}_d(\R)
\Bigg\}
$$
the \emph{statistical manifold} of centered Gaussians and simply write $\rho=\mathcal N(0,A)=\N(A)$.
Note that, for \emph{positive}-definite matrices $A\in\S^{++}_d(\R)$, the space of infinitesimal perturbations (the tangent plane $T_A\S^{++}_d(\R)$) is the whole space of symmetric matrices $U\in\S_d(\R)$.
\emph{Semi}-definite matrices are somehow degenerate (extremal) points of $\S^{++}_d(\R)$, but it is worth stressing that that \eqref{eq:def_Bures_static} makes sense even for \emph{semi}-definite $A_0,A_1\in\S^+_d(\R)$.

Going back to optimal transport, it is well-known \cite{takatsu2011wasserstein,modin} that the statistical manifold $\Stat_0$ is a totally geodesic submanifold in the Wasserstein space.
In other words, if $\rho_0=\N(A_0)$ and $\rho_1=\N(A_1)$ are Gaussians, then minimizers of the dynamical problem \eqref{eq:def_W2_dynamic} remain Gaussian, i-e $\rho_t=\mathcal N(0,A_t)$ for some suitable $A_t\in \S^{++}_d$.
Thus one expects that
$\Stat_0$, or rather $\S^{++}_d(\R)$, can be equipped with a well-chosen, finite-dimensional Riemannian metric such that the induced Riemannian distance coincides with the Wasserstein distance between corresponding Gaussians.
This is indeed the case \cite{takatsu2011wasserstein}, and the Riemannian metric is explicitly given at a point $\rho=\N(A)$ by
$$
g_A(U,V)=\tr(UAV)=\langle AU,V\rangle=\langle AV,U\rangle,
\qquad U,V\in T_A\S^{++}_d(\R)\cong \S_d(\R).
$$
Note that, if $g_{Eucl}(\xi,\zeta)=\langle \xi,\zeta\rangle=\tr (\xi \zeta^t)$ is the Euclidean scalar product, the corresponding Riesz isomorphism $U\mapsto \xi_U$ (at a point $A\in\S_d^{++}$) identifying $g_{Eucl}(\xi_U,\xi_V)=g_A(U,V)$ is given by
$
\xi_U=(AU)^{sym}=\frac{AU+UA}{2}
$.
(This is the exact equivalent of the elliptic correspondence \eqref{eq:Onsager_W} for optimal transport.)
As a consequence the Bures-Wasserstein distance can be computed \cite{takatsu2011wasserstein,MV} as
\begin{equation}
4\B^2(A_0,A_1)=
\min\limits_{A,U}\Bigg\{
\int_0^1 \langle A_tU_t,U_t\rangle\,\rd t
\qtext{s.t.}\frac {dA_t}{dt}=(A_tU_t)^{sym}\Bigg\},
\quad A_0,A_1\in\S^{++}_d
\label{eq:def_B_dynamic}
\end{equation}
Up to the scaling factor $4$ this an exact counterpart of \eqref{eq:def_W2_dynamic}, where $\frac {dA_t}{dt} =(A_tU_t)^{sym}$ plays the role of the continuity equation $\partial_t\rho_t+\dive(\rho_tv_t)=0$ and $\langle A_tU_t,U_t\rangle$ substitutes for the kinetic energy density $\int_\Omega\rho_t|v_t|^2$.

\section{The non-commutative Fisher-Rao space}
\label{sec:FR}

Given two probability densities $\rho_0,\rho_1\in\P(\D)$ on a domain $\D\subset \R^N$ (not to be confused with the previous $\Omega=\R^d$ for Gaussian optimal transport, codomain over which our matrices $A\in \S^{++}_d$ were built), the \emph{scalar} Fisher-Rao distance is classically defined as the Riemannian distance induced by the Fisher information metric
\begin{multline}
\FR^2(\rho_0,\rho_1)
=\min\limits_{\rho}\Bigg\{
\int_0^1 \int_\D \left|\frac{\partial\log\rho_t(x)}{\partial t}\right|^2  \rho_t(x)\,\rd x\,\rd t,
\quad\mbox{s.t. }\rho_t\in\P(\D)
\Bigg\}
\\
=\min\limits_{\rho,u}
\Bigg\{
\int_0^1 \int_\D \left|u_t(x)\right|^2  \rho_t(x)\,\rd x\,\rd t
\quad\mbox{s.t. }\partial_t\rho_t=\rho_t u_t
\mbox{ and }\rho_t\in\P(\D)
\Bigg\}.
\label{eq:FR_dynamic}
\end{multline}
Note that the unit-mass condition $\rho_t\in\P(\D)$ is enforced here as a hard constraint.
(Without this constraint, \eqref{eq:FR_dynamic} actually defines the \emph{Hellinger distance} $\mathsf H$ between arbitrary nonnegative measures \cite{MV}.)

The strong structural similarity between \eqref{eq:def_B_dynamic} \eqref{eq:FR_dynamic} suggest a natural extension of the latter scalar setting to \emph{matrix-valued} measures.
More precisely, the \emph{Fisher-Rao space} is the space of ($d$-dimensional) positive \emph{semi}-definite valued measures over $\D$ with unit mass
$$
\Pp(\D)=\Bigg\{ A\in\mathcal M(\D;\S^+_d(\R))\qtext{s.t.}\int_\D\tr A(x)\,\rd x=1\Bigg\}.
$$
This can be thought of as the space of noncommutative probability measures (see \cite{MV} for a short discussion on the connections with \emph{free probability} theory and $C^*$-algebras).
The \emph{matricial} Fisher-Rao distance is then defined as
\begin{multline}
 \FR^2(A_0,A_1)=\min\limits_{A,U}\Bigg\{\int_0^1\int_\D \langle A_t(x)U_t(x),U_t(x)\rangle\,\rd x\,\rd t
 \\
 \qtext{s.t.}
 \partial_tA_t=(A_tU_t)^{sym}
 \mbox{ in }(0,1)\times\D
 \qtext{and}
 A_t\in\Pp(\D)
 \Bigg\},
\label{eq:FR^2_dyn}
\end{multline}
which is indeed a higher-dimensional extension of \eqref{eq:FR_dynamic} and should be compared again with \eqref{eq:def_W2_dynamic}. 
This clearly suggests viewing $\Pp(\D)$ as a (formal, infinite-dimensional) Riemannian manifold with tangent space and norm
$$
T_A\Pp=\Bigg\{
\xi_U=(AU)^{sym}:
\quad
U\in L^2\left(A(x)\rd x;\,\S^+_d\right)
\mbox{ and }
\int_\D \langle A(x),U(x)\rangle\,\rd x =0
\Bigg\}
$$
\begin{equation}
\|\xi_U\|_{A}^2=\|U\|_{L^2_A}^2 =\int_\D \langle A(x)U(x),U(x)\rangle\,\rd x.
\label{eq:def_FR_speed}
\end{equation}
Note once again that we imposed the unit-mass condition $A_t\in\Pp$ as a hard constraint in \eqref{eq:FR^2_dyn}, which results in the zero-average condition $\int \langle A,U\rangle=0$ for tangent vectors $U$.
Removing the mass constraint leads to (an extension of) the Hellinger distance $\mathsf H$ between arbitrary PSD matrix-valued measures. (The Hellinger space $(\mathbb H^+,\mathsf H)$ is a metric cone over the ``unit-sphere'' $(\Pp,\FR)$ but we shall not discuss this rich geometry any further, see again \cite{MV}.)

The corresponding Riemannian gradients of internal-energy functionals $\mathcal F(A)=\int_\D F(A(x))\,\rd x$ can be explicited \cite{MV} as
\begin{equation}
\grad_{\FR}\mathcal F(A)=\left(A(x) F'(A(x))\right)^{sym}-\tr\left(\int_\D A(y)F'(A(y))\rd y\right)A(x).
\label{eq:grad_FR}
\end{equation}
Here $F'(A)$ stands for the usual first variation of the matricial function $F(A)$, computed with respect to the standard Euclidean (Frobenius) scalar product.
%
\section{Extended entropy, the heat flow, and the Schr\"odinger problem}
\label{sec:extended_entropy}
From now on we always endow $\Pp$ with the Riemannian metric given by \eqref{eq:def_FR_speed}.
Once this metric is fixed, the geometric Schr\"odinger problem \eqref{eq:geometric_schro} will be fully determined as soon as we choose a driving potential $V$ on $\Pp$.
One should ask now: 
\begin{center}
{\it
What is a good choice of ``the'' canonical entropy on the Fisher-Rao space?} 
\end{center}
A first natural guess would be based on the classical (negative) \emph{von~Neumann entropy} $S(A)=\tr (A\log A)$ from quantum statistical mechanics.
However, the functional $\mathcal S(A)=\int_\D S(A(x))\,\rd x$ lacks geodesic convexity with respect to our ambient Fisher-Rao metric \cite{MV}, and this makes it unfit for studying a suitable Schr\"odinger problem (see \cite{MVT} for the connection between the validity of the Gamma-convergence in the limit $\eps\to 0$ in \eqref{eq:geometric_schro} and the necessity of geodesic convexity of $V$ in a fairly general metric context).

It turns out that there is a second natural choice, dictated by our previous construction based on Gaussian optimal transport and the canonical Boltzmann entropy $H(\rho)=\int\rho\log\rho$ on $\P(\R^d)$.
Indeed, consider $A\in \S^{++}_d$ and write
$\rho_A=\N(A)$ for the corresponding Gaussian.
Assuming for simplicity that $\D$ has measure $|\D|=\frac 1d$ so that $\tr\int_\D \Id \rd x=1$, we choose as a reference measure the generalized uniform Lebesgue measure $\Id\in\Pp$ (more general reference measures can also be covered, see \cite{MV}).
An explicit computation then gives the  Boltzmann entropy of $\rho_A$ relatively to $\rho_I=\N(\Id)$ as
\begin{equation}
E(A)=H(\rho_A|\rho_I)
=\int_{\R^d}\frac{\rho_A(y)}{\rho_I(y)}\log\left(\frac{\rho_A(y)}{\rho_I(y)}\right)\rho_I(y)\,\rd y
=\frac 12\tr [A-\log A-\Id].
\label{eq:def_E}
\end{equation}
Note carefully that $E(A)=+\infty$ as soon as $A$ is only positive \emph{semi}-definite (due to $-\tr\log A=-\sum\log\lambda_i=+\infty$ if any of the eigenvalues $\lambda_i=0$), and that by convexity $E(A)\geq E(\Id)=0$ is minimal only when $A=\Id$. 
Our canonical definition of \emph{the} entropy on the Fisher-Rao space is then simply
\begin{equation}
\mathcal E(A)=\int_\D E(A(x))\,\rd x
= -\frac 12\tr \int_\D\log A(x)\,\rd x,
\qquad
\mbox{for }A=A(x)\in\Pp.
\label{eq:def_entropy_FR}
\end{equation}
(The terms $\int\tr A=1=\int\tr\Id$ cancel out in \eqref{eq:def_E} due to our mass normalization.)
Since the (Euclidean) first variation of \eqref{eq:def_E} is $E'(A)=\frac 12(\Id-A^{-1})$, \eqref{eq:grad_FR} shows that the Fisher-Rao gradient flow $\frac{dA_t}{dt}=-\grad_{\FR} \E(A_t)$ reads
\begin{multline}
\partial_t A_t  =-\left(A_t \frac 12[\Id-A_t^{-1}]\right)^{sym}
+\tr\left(\int_\D A_t \frac 12[\Id-A_t^{-1}]\rd y\right)A_t
\\
 = \frac 12(\Id-A_t)^{sym} + \frac 12\tr \left(\int _\D (A_t-\Id)\rd y\right) A_t = \frac 12(\Id-A_t).
 \label{eq:heat_ODE}
\end{multline}
This generalized ``heat flow'' is consistent with the observation that, whenever $\rho_0=\N(A_0)$ is Gaussian, the solution $\rho_t=\N(A_t)$ of the Fokker-Planck equation (the Wasserstein gradient flow of the relative entropy $\rho\mapsto H(\rho|\rho_I)$)
\begin{equation}
\partial_t\rho=\Delta\rho-\dive(\rho\nabla \log\rho_I)
\label{eq:Fokker-Planck}
\end{equation}
remains Gaussian with precisely $\frac{dA_t}{dt}=\frac 12(\Id-A_t)$.
Our extended Fisher information functional can then be defined as the dissipation rate of the entropy $\E$ along the ``heat flow'' $\frac{dA_t}{dt}=-\grad_{\FR} \E(A_t)$ given by \eqref{eq:heat_ODE}, namely
\begin{multline}
\mathcal F(A_t)=-\frac{d}{dt}\E(A_t)
=
\frac 12\frac{d}{dt}\tr\int_\D  \log A_t(x)\,\rd x
= \frac 12\tr \int_\D A_t^{-1}(x)\partial_tA_t(x)\,\rd x
\\
= \frac 12\tr \int_\D A_t^{-1}(x)\frac 12[\Id-A_t(x)]\,\rd x
=\frac 14\left(\tr\int_\D A_t^{-1}(x)\,\rd x -1\right).
\label{eq:def_extended_Fisher_info}
\end{multline}
Equivalently and consistently, $\mathcal F(A)=\|\grad_\FR \E(A)\|^2$.
Note that $\mathcal F(A_t)>F(\Id)=0$ and $\E(A_t)>\E(\Id)=0$ unless $A_t(x)\equiv\Id$, which is of course consistent with the expected long-time behavior $A_t\to \Id$ for \eqref{eq:heat_ODE} as $t\to\infty$ (or equivalently $\rho_t\to \rho_I=\N(\Id)$ for Gaussian solutions of the Fokker-Planck equation \eqref{eq:Fokker-Planck}).

With these explicit representations \eqref{eq:def_FR_speed}\eqref{eq:def_extended_Fisher_info} of the quantum (i-e matricial) Fisher-Rao metric and Fisher information, we can now make sense of the geometric Schr\"odinger problem \eqref{eq:geometric_schro} in the noncommutative Fisher-Rao space as
\begin{multline}
\min\limits_A
\Bigg\{
\int_0^1\int_\D \langle A_t(x)U_t(x),U_t(x)\rangle\,\rd x\rd t +\frac{\eps^2}4 \int_0^1 \left(\tr\int_\D A_t^{-1}(x)\,\rd x -1\right)
\\
\mbox{s.t. }\qquad \partial_tA_t=(A_t U_t)^{sym}
\qtext{and}A_t\in\Pp\mbox{ for all }t\in[0,1]
\Bigg\},
\label{eq:Schro_FR}
\end{multline}
with fixed endpoints $A_0,A_t\in \Pp(\D)$.
Let us mention at this stage that a different entropic regularization of Gaussian optimal transport was investigated in \cite{JMPC} in a static framework and for a much simpler setting, namely when $\D=\{x\}$ is a single point

As could be expected from the above geometric machinery, we have now
\begin{theo*}[\cite{MV}]
 In the limit $\eps\to 0$ and for fixed endpoints $A_0,A_1\in\Pp$, the $\eps$-functional in \eqref{eq:Schro_FR} Gamma-converges towards the kinetic functional in \eqref{eq:FR^2_dyn} for the uniform convergence on $C([0,1];\Pp)$.
\end{theo*}
We omit the details for the sake of brevity, but let us point out that our proof leverages fully explicit properties of the geometric heat flow \eqref{eq:heat_ODE}.
A particular byproduct of our analysis is the $\frac 12$-geodesic convexity of the entropy \eqref{eq:def_entropy_FR} in the Fisher-Rao geometry, which was already established in \cite{modin} by formal Riemannian computations.
The key argument builds up on a Lagrangian construction originally due to A. Baradat in \cite{BM}.
Finally, we recently extended the result to arbitrary metric spaces and general entropy functionals \cite{MVT} (provided a suitable heat flow is available), and we showed moreover that geodesic convexity is actually necessary (and almost sufficient) for the $\Gamma$-convergence.
%
%
%
%

\end{document}